\documentclass{amsart}

\theoremstyle{definition}

\theoremstyle{remark}

\reversemarginpar

\begin{document}

\title[Digamma function]{The integrals in Gradshteyn and Ryzhik. \\
Part 10: The digamma function}

\subjclass[2000]{Primary 33}


\author[L. Medina]{Luis A. Medina}
\address{Department of Mathematics,
Tulane University, New Orleans, LA 70118}
\email{lmedina@math.tulane.edu}

\author[V. Moll]{Victor H. Moll}
\address{Department of Mathematics,
Tulane University, New Orleans, LA 70118}
\email{vhm@math.tulane.edu}

\thanks{The second author wishes to acknowledge the partial support of  
NSF-DMS 0409968. The first author was partially supported as a graduate
student by the same grant.}

\begin{abstract}
The table of Gradshteyn and Rhyzik contains some  integrals  that can be 
expressed in terms of the digamma function $\psi(x) = \frac{d}{dx} \log 
\Gamma(x)$. In this note we present some of these evaluations.
\end{abstract}

\maketitle

\newcommand{\nn}{\nonumber}
\newcommand{\ba}{\begin{eqnarray}}
\newcommand{\ea}{\end{eqnarray}}
\newcommand{\ift}{\int_{0}^{\infty}}
\newcommand{\ione}{\int_{0}^{1}}
\newcommand{\ifft}{\int_{- \infty}^{\infty}}
\newcommand{\no}{\noindent}
\newcommand{\realpart}{\mathop{\rm Re}\nolimits}
\newcommand{\imagpart}{\mathop{\rm Im}\nolimits}

\newtheorem{Definition}{\bf Definition}[section]
\newtheorem{Thm}[Definition]{\bf Theorem} 
\newtheorem{Example}[Definition]{\bf Example} 
\newtheorem{Lem}[Definition]{\bf Lemma} 
\newtheorem{Note}[Definition]{\bf Note} 
\newtheorem{Cor}[Definition]{\bf Corollary} 
\newtheorem{Prop}[Definition]{\bf Proposition} 
\newtheorem{Problem}[Definition]{\bf Problem} 
\numberwithin{equation}{section}

\section{Introduction} \label{intro} 
\setcounter{equation}{0}

The table of integrals \cite{gr} contains a large variety of definite 
integrals that involve the {\em digamma} function 
\begin{equation}
\psi(x) = \frac{d}{dx} \log \Gamma(x) = \frac{\Gamma'(x)}{\Gamma(x)}.
\end{equation}
\noindent
Here $\Gamma(x)$ is the gamma function defined by 
\begin{equation}
\Gamma(x) = \ift t^{x-1}e^{-t} \, dt. 
\end{equation}

Many of the analytic properties can be derived from those of $\Gamma(x)$. 
The next theorem represents a collection of the important properties of 
$\Gamma(x)$ that are used in the current paper. The reader will find in
\cite{irrbook} detailed proofs. 

\begin{Thm}
The gamma function satisfies: \\

\noindent
a) the functional equation
\begin{equation}
\Gamma(x+1) = x \Gamma(x). 
\end{equation}

\medskip

\noindent
b) For $n \in \mathbb{N}$, the interpolation formula 
$\Gamma(n) = (n-1)!$.

\medskip

\noindent
c) The {\em Euler constant} $\gamma$, defined by 
\begin{equation}
\gamma = \lim\limits_{n \to \infty} \sum_{k=1}^{n} \frac{1}{k} - \ln n,
\end{equation}
\noindent
is also given by $\gamma = - \Gamma'(1)$.  This appears as the special case
$a = 1$ of formula $\mathbf{4.331.1}$: 
\begin{equation}
\ift e^{-ax} \ln x \, dx = -\frac{\gamma + \ln a}{a}. 
\end{equation}
\noindent
This was established in \cite{moll-gr4}. 
The change of variables $t = ax$ shows that the case $a=1$ is equivalent to 
the general case. This is an instance of a {\em fake parameter}. 

\medskip

\noindent
d) The infinite product representation
\begin{equation}
\Gamma(x) = \frac{e^{-\gamma x}}{x} \prod_{k=1}^{\infty} 
\left[ \left( 1 + \frac{x}{k} \right)^{-1} e^{x/k} \right]
\label{gamma-product}
\end{equation}
\noindent
is valid for $x \in \mathbb{C}$ away from the poles at $x=0, \, -1, \, -2, 
\ldots$

\medskip
\noindent
e) For $n \in \mathbb{N}$ we have
\begin{equation}
\Gamma \left( n + \tfrac{1}{2} \right) = \frac{(2n)!}{2^{2n} \, n!} \sqrt{\pi}
\label{gamma-half}
\end{equation}
\noindent
and
\begin{equation}
\Gamma  \left(\tfrac{1}{2}-n \right) = (-1)^{n} \frac{2^{2n} \, n!}{(2n)!} 
\, \sqrt{\pi}.
\end{equation}

\medskip
\noindent
f) For $x \in \mathbb{C}, \, x \not \in \mathbb{Z}$ we have the 
reflection rule
\begin{equation}
\Gamma(x) \Gamma(1-x) = \frac{\pi}{\sin \pi x}. 
\label{reflec-gamma}
\end{equation}
\end{Thm}

\medskip

Several properties of the digamma function $\psi(x)$ follow directly from the
gamma function. 

\begin{Thm}
The digamma function $\psi(x)$ satisfies 

\noindent
a) the functional equation 
\begin{equation}
\psi(x+1)  = \psi(x) + \frac{1}{x}. 
\label{psi1}
\end{equation}

\medskip

\noindent
b) For $n \in \mathbb{N}$, we have 
\begin{equation}
\psi(n) = - \gamma + \sum_{k=1}^{n-1} \frac{1}{k}. 
\label{psi-n}
\end{equation}
\noindent
In particular, $\psi(1) = - \gamma$. 

\medskip

\noindent
c) For $x \in \mathbb{C}$ away from $x=0, \, -1, \, -2, \, \cdots$ we have
\begin{eqnarray}
\psi(x) & = & -\gamma -\frac{1}{x} + x \sum_{k=1}^{\infty} \frac{1}{k(x+k)}
\label{psi-series} \\
& = & -\gamma - \sum_{k=0}^{\infty} \left( \frac{1}{x+k} - \frac{1}{k+1} 
\right) \nonumber
\end{eqnarray}

\medskip

\noindent
d) The derivative of $\psi$ is given by
\begin{equation}
\psi'(x) = \sum_{k=0}^{\infty} \frac{1}{(x+k)^{2}}. 
\end{equation}
\noindent
In particular, $\psi'(1) = \pi^{2}/6$. 

\medskip

\noindent
e) For $n \in \mathbb{N}$ we have
\begin{equation}
\psi( \tfrac{1}{2} \pm n  ) = - \gamma - 2 \ln 2 + 2 \sum_{k=1}^{n} 
\frac{1}{2k-1}.
\label{psi-half}
\end{equation}
\noindent
In particular, 
\begin{equation}
\psi( \tfrac{1}{2} ) = -\gamma - 2 \ln 2.
\end{equation}

\medskip
f) For $x \in \mathbb{C}, \, x \not \in \mathbb{Z}$ we have the reflection 
rule
\begin{equation}
\psi(1-x)  = \psi(x) + \pi \, \text{cot } \pi x. 
\label{psi2}
\end{equation}
\end{Thm}

\section{A first integral representation} \label{sec-first} 
\setcounter{equation}{0}

In this section we establish the integral evaluation $\mathbf{3.429}$. 
Severeal direct consequences of this formulas are also described. 

\begin{Prop}
Assume $a>0$. Then 
\begin{equation}
\ift \left[ e^{-x} - (1+x)^{-a} \right] \frac{dx}{x} = \psi(a).
\label{psi-1}
\end{equation}
\end{Prop}
\begin{proof}
We begin with the double integral
\begin{equation}
\ift \int_{1}^{s} e^{-tz} dt \, dz = 
\ift \frac{e^{-z} - e^{-sz}}{z} dz. 
\end{equation}
\noindent
On the other hand,
\begin{equation}
\int_{1}^{s} \ift e^{-tz} dz \, dt =  \int_{1}^{s} \frac{dt}{t} = \ln s.
\end{equation}
\noindent
We conclude that
\begin{equation}
\ift \frac{e^{-z} - e^{-sz}}{z} \, dz = \ln s.
\end{equation}
\noindent
This evaluation is equivalent to:
\begin{equation}
\ift \frac{e^{-ax} - e^{-bx}}{x} \, dx = \ln \frac{b}{a},
\label{form25}
\end{equation}
\noindent
that appears as formula $\mathbf{3.434.2}$ in \cite{gr}. The reader 
will find a proof in \cite{moll-gr4}.

We now establish the result: start with 
\begin{eqnarray}
\Gamma'(a)  & = & \ift e^{-s} s^{a-1} \ln s \, ds \nonumber \\
 & = & \ift e^{-s} s^{a-1} \ift \frac{e^{-z} - e^{-zs}}{z} \, dz \, ds 
\nonumber \\
& = & \ift \left( e^{-z} \ift s^{a-1} e^{-s} \, ds - 
\ift s^{a-1} e^{-s(1+z)} \, ds \right) \, \frac{dz}{z}. \nonumber
\end{eqnarray}
\noindent
This formula can be rewritten as 
\begin{equation}
\Gamma'(a) = \Gamma(a) \, \ift \left( e^{-z} - (1+z)^{-a} \right) 
\frac{dz}{z}. \nonumber
\end{equation}
\noindent
This establishes (\ref{psi-1}).  
\end{proof}

\medskip

\begin{Example}
\label{example1}
The special case $a=1$ yields
\begin{equation}
\ift \left( e^{-x} - \frac{1}{1+x} \right) \, \frac{dx}{x} = - \gamma.
\label{34353}
\end{equation}
\noindent
This appears as $\mathbf{3.435.3}$. 
\end{Example}

\begin{Example}
The change of variables $w = - \ln x$ gives the value of $\mathbf{4.275.2}$:
\begin{equation}
\int_{0}^{1} \left[ x  - \left( \frac{1}{1- \ln x} \right)^{q} \right] 
\frac{dx}{x \, \ln x} = 
- \int_{0}^{\infty} \left[ e^{-w} - (1+w)^{-q} \right] \frac{dw}{w} = 
- \psi(q).
\end{equation}
\end{Example}

\medskip

\begin{Example}
The change of  variables $ t = 1/(x+1)$ in (\ref{psi-1}) yields 
$\mathbf{3.471.14}$:
\begin{equation}
\int_{0}^{1} \frac{e^{(1-1/t)} - t^{a}}{t(1-t)} \, dt = \psi(a)
\end{equation}
\end{Example}

\medskip

\begin{Example}
The result of Example \ref{example1} can be used to prove 
$\mathbf{3.435.4}$:
\begin{equation}
\ift \left( e^{-bx} - \frac{1}{1+ax} \right) \frac{dx}{x} = 
\ln \frac{a}{b} - \gamma. 
\end{equation}
\noindent
Indeed, the change of variables $t = bx$ yields
from (\ref{example1}) the identity

\begin{eqnarray}
\ift \left( e^{-bx} - \frac{1}{1+ax} \right) \frac{dx}{x} & = &  
\ift \left( e^{-t} - \frac{1}{1+at/b} \right) \frac{dt}{t} \nonumber \\
& = & \ift \frac{e^{-t} - e^{-at/b}}{t} \, dt  + 
\ift \left( e^{-at/b} - \frac{1}{1+at/b} \right) \frac{dt}{t}. \nonumber
\end{eqnarray}
\noindent
Formula (\ref{form25}) shows the first integral is $\ln \frac{a}{b}$ and 
the value of the second one comes from (\ref{example1}). 
\end{Example}

\begin{Example}
The evaluation $\mathbf{3.476.2}$:
\begin{equation}
\int_{0}^{\infty} \left( e^{-x^{p}} - e^{-x^{q}} \right) 
\frac{dx}{x} = \frac{p-q}{pq} \gamma
\label{34762}
\end{equation}
\noindent
comes directly from (\ref{psi-1}). Indeed, the change of variables 
$u = x^{p}$ yields
\begin{equation}
I := \int_{0}^{\infty} \left(e^{-x^{p}} - e^{-x^{q}} \right)
\frac{dx}{x} = 
\frac{1}{p} \int_{0}^{\infty} \left( e^{-u} - e^{-u^{q/p}} \right) \frac{du}{u}.
\nonumber
\end{equation}
\noindent
Now write
\begin{equation}
I = \frac{1}{p} \int_{0}^{\infty} \left( e^{-u} - \frac{1}{1+u} \right) 
\frac{du}{u} + \frac{1}{p} \int_{0}^{\infty} \left( \frac{1}{1+u} - 
e^{-u^{q/p}} \right) \frac{du}{u}. \nonumber
\end{equation}
\noindent
The first integral is $-\gamma$ by (\ref{34353})  and the change of variables
$v = u^{q/p}$ gives
\begin{eqnarray}
I & = & -\frac{\gamma}{p} + \frac{1}{q} 
\int_{0}^{\infty} \left( \frac{1}{1+ v^{p/q}} - e^{-v} \right) \frac{dv}{v}
\nonumber \\
& = & -\frac{\gamma}{p} 
+ \frac{1}{q} \int_{0}^{\infty} \left( \frac{1}{1+v} - e^{-v} \right) 
\frac{dv}{v} + 
\frac{1}{q} \int_{0}^{\infty} \frac{v - v^{p/q}}{v(1+v)(1+v^{p/q})} \, dv.
\nonumber
\end{eqnarray}
\noindent
Split the last integral from $[0,1]$ to $[1, \infty)$ and use the change 
of variables $x \mapsto 1/x$ in the second part to check that the whole
integral vanishes. Formula (\ref{34762}) has been established.
\end{Example}

\medskip

\begin{Example}
Formula $\mathbf{3.463}$:
\begin{equation}
\int_{0}^{\infty} \left( e^{-x^{2}} - e^{-x} \right) \frac{dx}{x} = 
\frac{\gamma}{2}
\end{equation}
\noindent
corresponds to the choice $p=2$ and $q=1$ in (\ref{34762}).
\end{Example}

\medskip

\begin{Example}
Formula $\mathbf{3.469.2}$:
\begin{equation}
\int_{0}^{\infty} \left( e^{-x^{4}} - e^{-x} \right) \frac{dx}{x} = 
\frac{3 \gamma}{4}
\end{equation}
\noindent
corresponds to the choice $p=4$ and $q=1$ in (\ref{34762}).
\end{Example}

\medskip

\begin{Example}
Formula $\mathbf{3.469.3}$:
\begin{equation}
\int_{0}^{\infty} \left( e^{-x^{4}} - e^{-x^{2}} \right) \frac{dx}{x} = 
\frac{\gamma}{4}
\end{equation}
\noindent
corresponds to the choice $p=4$ and $q=2$ in (\ref{34762}).
\end{Example}

\medskip

\begin{Example}
Formula $\mathbf{3.475.3}$:
\begin{equation}
\int_{0}^{\infty} \left( e^{-x^{2^{n}}} - e^{-x} \right) \frac{dx}{x} = 
(1- 2^{-n}) \gamma
\end{equation}
\noindent
corresponds to the choice $p=2^{n}$ and $q=1$ in (\ref{34762}).
\end{Example}

\medskip

The case $p=q$ in (\ref{34762}) is now modified to include a paramter. 

\begin{Prop}
Let $a, \, b, \, p \in \mathbb{R}^{+}$. Then $\mathbf{3.476.1}$ in \cite{gr}
states that
\begin{equation}
\int_{0}^{\infty} \left[ e^{-ax^{p}} - e^{-bx^{p}} \right] \frac{dx}{x} =
\frac{\ln b - \ln a}{p}.
\end{equation}
\end{Prop}
\begin{proof}
The change of variables $t = ax^{p}$ gives
\begin{equation}
\int_{0}^{\infty} \left[ e^{-ax^{p}} - e^{-bx^{p}} \right] \frac{dx}{x} =
\frac{1}{p} \int_{0}^{\infty} \left( e^{-t} - e^{-bt/a} \right) \frac{dt}{t}.
\nonumber
\end{equation}
\noindent
Introduce the term $1/(1+t)$ to obtain
\begin{eqnarray}
I & = & \frac{1}{p} \int_{0}^{\infty} \left( e^{-t} - \frac{1}{1+t} \right) 
\frac{dt}{t} - 
\frac{1}{p} \int_{0}^{\infty} \left( e^{-bt/a} - \frac{1}{1+t} \right) 
\frac{dt}{t} \nonumber \\
 & = & - \frac{\gamma}{p} - \frac{1}{p} 
\int_{0}^{\infty} \left( e^{-s} - \frac{b}{b+ as} \right) \frac{ds}{s}.
\nonumber
\end{eqnarray}
\noindent
Adding and subtracting the term $1/(1+s)$  produces
\begin{equation}
I = \frac{1}{p} \int_{0}^{\infty} 
\left( \frac{b}{b + a s} - \frac{1}{1+s} \right) \frac{ds}{s}.
\end{equation}
\noindent
The final result now comes from evaluating the last integral. 
\end{proof}

\medskip

We now present another integral representation of the digamma function. \\

\begin{Prop}
The digamma function is given by 
\begin{equation}
\psi(a) = \ift \left( \frac{e^{-x}}{x} - \frac{e^{-ax}}{1-e^{-x}} \right) \, dx.
\label{psi-2}
\end{equation}
\noindent
This expression appears as $\mathbf{3.427.1}$ in \cite{gr}. 
\end{Prop}
\begin{proof}
The representation (\ref{psi-1}) is written  as
\begin{equation}
\psi(a) = \lim\limits_{\delta \to 0} \int_{\delta}^{\infty} \frac{e^{-z}}{z}
\, dz - \int_{\delta}^{\infty} \frac{dz}{z(1+z)^{a}},
\end{equation}
\noindent
to avoid the singularity at $z=0$. The change of variables $z=e^{t}-1$ in 
the second integral gives
\begin{equation}
\psi(a) = \lim\limits_{\delta \to 0} \int_{\delta}^{\infty} \frac{e^{-z}}{z}
\, dz - \int_{\ln(1+\delta)}^{\infty} \frac{e^{-at} \, dt}{1-e^{-t}}.
\end{equation}
\noindent
Now observe that
\begin{equation}
\left|\int_{\delta}^{\ln(1+\delta)} \frac{e^{-t}}{t} \, dt \right|
\leq \int_{\ln(1+\delta)}^{\delta} \frac{dt}{t} \to 0, 
\end{equation}
\noindent
as $\delta \to 0$. This completes the proof. 
\end{proof}

\begin{Example}
The special case $a=1$ in (\ref{psi-2}) gives $\mathbf{3.427.2}$:
\begin{equation}
\ift \left( \frac{1}{1-e^{-x}} - \frac{1}{x} \right)e^{-x} \, dx = \gamma.
\end{equation}
\end{Example}

\medskip

\begin{Example}
The change  of variables $t = e^{-x}$ in (\ref{psi-2}) produces
$\mathbf{4.281.4}$:
\begin{equation}
\int_{0}^{1} \left( \frac{1}{\ln t} + \frac{t^{a-1}}{1-t} \right) \, dt 
= - \psi(a). 
\label{psi-3}
\end{equation}
\end{Example}

\begin{Example}
The special case $a=1$ in (\ref{psi-3}) yields $\mathbf{4.281.1}$:
\begin{equation}
\int_{0}^{1} \left( \frac{1}{\ln t} + \frac{1}{1-t} \right) \, dt = \gamma. 
\end{equation}
\end{Example}

\begin{Prop}
Let $p, \, q \in \mathbb{R}$. Then 
\begin{equation}
\int_{0}^{1} \left( \frac{x^{p-1}}{\ln x } + \frac{x^{q-1}}{1-x} \right) \, dx 
= \ln p - \psi(q).
\label{form42815}
\end{equation}
\noindent
This appears as $\mathbf{4.281.5}$ in \cite{ gr}.
\end{Prop}
\begin{proof}
Write 
\begin{equation}
\int_{0}^{1} \left( \frac{x^{p-1}}{\ln x } + \frac{x^{q-1}}{1-x} \right) \, dx 
= \int_{0}^{1} \left( \frac{1}{\ln x } + \frac{x^{q-1}}{1-x} \right) \, dx 
+ \int_{0}^{1} \frac{x^{p-1}-1}{\ln x } \, dx.
\end{equation}
\noindent
The first integral is $- \psi(q)$ from (\ref{psi-3}) and to evaluate the 
second one, differentiate with respect to $p$, to produce
\begin{equation}
\frac{d}{dp} \int_{0}^{1} \frac{x^{p-1}-1}{\ln x } \, dx = 
\int_{0}^{1}x^{p-1} \, dx = \frac{1}{p}.
\end{equation}
\noindent
The value at $p=1$ shows that the constant of integration vanishes. The 
formula (\ref{form42815}) has been established. 
\end{proof}

\section{The difference of values of the digamma function} \label{sec-simple} 
\setcounter{equation}{0}

In this section we establish an integral representation for the difference of 
values of the digamma function. The expression appears as 
$\mathbf{3.231.5}$ in \cite{gr}.

\begin{Prop}
Let $p, \, q \in \mathbb{R}$. Then 
\begin{equation}
\int_{0}^{1} \frac{x^{p-1} - x^{q-1}}{1-x} \, dx = \psi(q) - \psi(p).
\label{3.231.5}
\end{equation}
\end{Prop}
\begin{proof}
Consider first 
\begin{equation}
I( \epsilon) = 
\int_{0}^{1} x^{p-1}(1-x)^{\epsilon-1} \, dx - 
\int_{0}^{1} x^{q-1}(1-x)^{\epsilon-1} \, dx,
\end{equation}
\noindent
that avoids the apparent 
singularity at $x=1$. The integral $I(\epsilon)$ can be 
expressed in terms of the beta function 
\begin{equation}
B(a,b) = \int_{0}^{1} x^{a-1}(1-x)^{b-1} \, dx
\end{equation}
\noindent
as $I( \epsilon) =  B(p, \epsilon) - B(q, \epsilon)$,  and using the relation
\begin{equation}
B(a,b) = \frac{\Gamma(a) \, \Gamma(b)}{\Gamma(a+b)} 
\end{equation}
\noindent 
we obtain 
\begin{equation}
I( \epsilon) =  \Gamma(\epsilon) \left( \frac{\Gamma(p)}{\Gamma(p+ \epsilon)}
- \frac{\Gamma(q)}{\Gamma(q + \epsilon)} \right). 
\end{equation}
\noindent
Now use $\Gamma(1+ \epsilon) = \epsilon \Gamma(\epsilon)$ to write
\begin{equation}
I( \epsilon) = \Gamma(1+ \epsilon) 
\left( \frac{\Gamma(p)-\Gamma(p + \epsilon)}{\epsilon} \, 
\frac{1}{\Gamma(p + \epsilon)} - 
 \frac{\Gamma(q)-\Gamma(q + \epsilon)}{\epsilon} \, 
\frac{1}{\Gamma(q + \epsilon)} \right), 
\end{equation}
\noindent
and obtain (\ref{3.231.5}) by letting $\epsilon \to 0$. 
\end{proof}

\medskip

\begin{Example}
The special value $\psi(1) = - \gamma$ produces 
\begin{equation}
\int_{0}^{1} \frac{1-x^{q-1}}{1-x} \, dx =  \gamma + \psi(q). 
\end{equation}
\noindent
This appears as $\mathbf{3.265}$ in \cite{gr}.
\end{Example}

\begin{Example}
A second special value appears in $\mathbf{3.268.2}$:
\begin{equation}
\int_{0}^{1} \frac{1-x^{a}}{1-x} \, x^{b-1} \, dx = 
\psi(a+b) - \psi(b).
\end{equation}
\noindent
It is obtained from (\ref{3.231.5}) by choosing $p = b$ and $q = a+b$. 
\end{Example}

\begin{Example}
Now let $q = 1-p$ in (\ref{3.231.5}) to produce
\begin{equation}
\int_{0}^{1} \frac{x^{p-1} - x^{-p}}{1-x} \, dx = \psi(1-p) - \psi(p) = 
\pi \, \text{cot }\pi p.
\label{3.231.1}
\end{equation}
\noindent
This appears as $\mathbf{3.231.1}$ in \cite{gr}.  
\end{Example}

\begin{Example}
The special 
case $p = a+1$ and $q= 1-a$ produces
\begin{equation}
\int_{0}^{1} \frac{x^{a} - x^{-a}}{1-x} \, dx = \psi(1-a) - \psi(1+a) =
\pi \, \text{cot }\pi a - \frac{1}{a}, 
\end{equation}
\noindent
where we have used (\ref{psi1}) and (\ref{psi2}) to simplify the result. 
This is $\mathbf{3.231.3}$ in \cite{gr}. 
\end{Example}

\begin{Example}
The change of variables $x = t^{a}$ in (\ref{3.231.5}) produces 
\begin{equation}
\int_{0}^{1} \frac{t^{ap-1}-t^{aq-1}}{1-t^{a}} \, dt = 
\frac{\psi(q) - \psi(p)}{a}. 
\label{scaled}
\end{equation}
\noindent
Now let $p=1, a = \nu$ and $q = \frac{\mu}{\nu}$  and the replace $\mu$ by $p$
and $\nu$ by $q$ to 
obtain $\mathbf{3.244.3}$ in 
\cite{gr}:
\begin{equation}
\int_{0}^{1} \frac{t^{q-1}-t^{p-1}}{1-t^{q}} \, dt = \frac{1}{q} 
\left( \gamma + \psi \left( \frac{p}{q} \right) \right). 
\label{32443}
\end{equation}
\end{Example}

\begin{Example}
The special case $p = b/a$ and $q = 1 - b/a$ in (\ref{scaled}) produces 
\begin{equation}
\int_{0}^{1} \frac{x^{b-1} - x^{a - b -1}}{1-x^{a}} \, dx = 
\frac{1}{a} \left( \psi( 1 - b/a) - \psi(b/a) \right). 
\end{equation}
\noindent
The result is  now simplified using (\ref{psi2}) to produce 
\begin{equation}
\int_{0}^{1} \frac{x^{b-1} - x^{a - b -1}}{1-x^{a}} \, dx = 
\frac{\pi}{a} \text{ cot } \frac{\pi b}{a}.
\end{equation}
\noindent
This is $\mathbf{3.244.2}$ in \cite{gr}. 
\end{Example}

\begin{Example}
The special case $a=2$ in (\ref{scaled}) yields
\begin{equation}
\int_{0}^{1} \frac{t^{2 \mu -1} - t^{2 \nu -1}}{1-t^{2}} \, dt = 
\frac{1}{2} \left( \psi(\nu) - \psi(\mu) \right). 
\label{scaled11}
\end{equation}
\noindent
The choice $\mu = 1 + p/2$ and $\nu = 1 -  p/2$:
\begin{equation}
\int_{0}^{1} \frac{x^{p}-x^{-p}}{1-x^{2}} \, x \, dx = 
\frac{1}{2} \left( \psi( 1 + p/2) - \psi(1 - p/2) \right). 
\end{equation}
\noindent
The identities $\psi(x+1) = \psi(x) + 1/x$ and 
$\psi(1-x) - \psi(x) = \pi \, \cot \pi x$ produce
\begin{equation}
\int_{0}^{1} \frac{x^{p}-x^{-p}}{1-x^{2}} \, x \, dx = 
\frac{\pi}{2} \cot \left( \frac{p \pi}{2} \right) - \frac{1}{p}. 
\end{equation}
\noindent
This appears as $\mathbf{3.269.1}$ in \cite{gr}.
\end{Example}

\begin{Example}
The choice $\mu = \frac{a+1}{2}$ and $\nu = \frac{b+1}{2}$ in 
(\ref{scaled11}) gives $\mathbf{3.269.3}$:
\begin{equation}
\int_{0}^{1} \frac{x^{a}-x^{b}}{1-x^{2}} \, dx = \frac{1}{2} 
\left( \psi \left( \frac{b+1}{2} \right) - \psi \left( \frac{a+1}{2} \right) \right).
\end{equation}
\end{Example}

\section{Integrals over a half-line} \label{sec-half} 
\setcounter{equation}{0}

In this section we consider integrals over the half-line $[0, \infty)$ that
can be evaluated in terms of the digamma function. 

\begin {Prop}
Let $p, \, q \in \mathbb{R}$. Then
\begin{equation}
\ift \left( \frac{t^{p}}{(1+t)^{p}} - \frac{t^{q}}{(1+t)^{q}} \right) \, 
\frac{dt}{t} = \psi(q) - \psi(p). 
\end{equation}
\noindent
This is $\mathbf{3.219}$ in \cite{gr}. Also
\begin{equation}
\ift \left( \frac{1}{(1+t)^{p}} - \frac{1}{(1+t)^{q}} \right) \, 
\frac{dt}{t} = \psi(q) - \psi(p). 
\label{32315}
\end{equation}
\end{Prop}
\begin{proof}
Let $t = x/(1-x)$ in (\ref{3.231.5}).  The second form comes from the 
first by the change of variables $x \mapsto 1/x$.
\end{proof}

\begin{Example}
The special case $p=1$  yields
\begin{equation}
\ift \left( \frac{1}{1+t} - \frac{1}{(1+t)^{q}} \right) \, 
\frac{dt}{t} = \psi(q)  + \gamma. 
\label{3.233}
\end{equation}
\noindent
This appears as $\mathbf{3.233}$ in \cite{gr}.
\end{Example}

\begin{Example}
The evaluation of $\mathbf{3.235}$:
\begin{equation}
\ift \frac{(1+x)^{a} - 1}{(1+x)^{b}} \frac{dx}{x} = \psi(b) - \psi(b-a)
\end{equation}
\noindent
can be established directly from (\ref{3.233}). Simply write
\begin{equation}
\ift \frac{(1+x)^{a} - 1}{(1+x)^{b}} \frac{dx}{x} = 
\ift \left( \frac{1}{1+x} - \frac{1}{(1+x)^{b}} \right)  \frac{dx}{x} - 
\ift \left( \frac{1}{1+x} - \frac{1}{(1+x)^{b-a}} \right)  \frac{dx}{x}
\nonumber
\end{equation}
\noindent
to obtain the result. 
\end{Example}

Some examples of integrals over $[0,\infty)$ can be reduced to a pair 
of integrals over $[0,1]$. 

\begin{Prop}
The formula $\mathbf{3.231.6}$ of \cite{gr} states that
\begin{equation}
\ift \frac{x^{p-1} -x^{q-1}}{1-x} \, dx = \pi \left( \cot \pi p - 
\cot \pi q \right).
\end{equation}
\end{Prop}
\begin{proof}
To evaluate this, make the change of variables $t = 1/x$ in the part over 
$[1, \infty)$ to produce
\begin{equation}
\ift \frac{x^{p-1} -x^{q-1}}{1-x} \, dx =
\int_{0}^{1}  \frac{x^{p-1} -x^{q-1}}{1-x} \, dx -
\int_{0}^{1}  \frac{x^{-p} -x^{-q}}{1-x} \, dx. 
\end{equation}
\noindent
Now use the result (\ref{3.231.5}) to write
\begin{equation}
\ift \frac{x^{p-1} -x^{q-1}}{1-x} \, dx =
\psi(q)-\psi(p) - \left[ \psi(1-q) - \psi(1-p) \right]. 
\end{equation}
\noindent
The relation $\psi(x)-\psi(1-x) = - \pi \, \cot(\pi q)$ yields the result. 
\end{proof}

\section{An exponential scale} \label{sec-expon} 
\setcounter{equation}{0}

In this section we present the evaluation of certain definite integrals 
involving the exponential function. These are integrals that can be 
evaluated in terms of the digamma function of the parameters involved. 

\begin{Example}
The simplest one is $\mathbf{3.317.2}$:
\begin{equation}
\int_{-\infty}^{\infty} \left( \frac{1}{(1+e^{-x})^{p}} - 
\frac{1}{(1+e^{-x})^{q}} \right) \, dx = \psi(q) - \psi(p)
\label{33172}
\end{equation}
\noindent
that comes from (\ref{32315}) via the change of variables $x \mapsto e^{-x}$. 
\end{Example}

\begin{Example}
The special case $p=1$ and $\psi(1) = - \gamma$ produces $\mathbf{3.317.1}$:
\begin{equation}
\int_{-\infty}^{\infty} \left( \frac{1}{1+e^{-x}} - 
\frac{1}{(1+e^{-x})^{q}} \right) \, dx = \psi(q) + \gamma
\end{equation}
\end{Example}

\begin{Example}
The evaluation of $\mathbf{3.316}$:
\begin{equation}
\int_{-\infty}^{\infty} \frac{(1+e^{-x})^{p} -1}{(1+e^{-x})^{q}} \, dx = 
\psi(q) - \psi(q-p)
\end{equation}
\noindent
comes directly from (\ref{33172}).
\end{Example}

\medskip

\begin{Prop}
Let $p, \, q \in \mathbb{R}$. Then 
\begin{equation}
\ift \frac{e^{-pt} - e^{-qt}}{1-e^{-t}} \, dt = \psi(q) - \psi(p),
\label{psiaa}
\end{equation}
\noindent
This appears as $\mathbf{3.311.7}$ in \cite{gr}. 
\end{Prop}
\begin{proof}
Make the change of variables $x = e^{-t}$ in (\ref{3.231.5}).
\end{proof}

\begin{Example}
The evaluation (\ref{psiaa}) can also be written as 
\begin{equation}
\ift \frac{e^{t(1-p)} - e^{t(1-q)}}{e^{t}-1} \, dt = \psi(q) - \psi(p),
\end{equation}
\end{Example}

\begin{Example}
The special case $p=1$ and $q= 1 - \nu$ is 
\begin{equation}
\ift \frac{1 - e^{\nu t}}{e^{t}-1} \, dt = \psi(1- \nu) - \psi(1),
\end{equation}
\noindent
and using $\psi(1) = -\gamma$ and $\psi(1-\nu) = \psi(\nu) + \pi \, 
\text{ cot } \pi \nu$, yields the form 
\begin{equation}
\ift \frac{1 - e^{\nu t}}{e^{t}-1} \, dt = \psi(\nu) + \gamma + 
\pi \, \text{ cot }\pi \nu,
\end{equation}
\noindent
as it appears in $\mathbf{3.311.5}$. 
\end{Example}

\begin{Example}
 Another special case of (\ref{psiaa}) is 
$\mathbf{3.311.6}$, that corresponds to $p=1$: 
\begin{equation}
\ift \frac{e^{-t} - e^{-qt}}{1 - e^{-t}} \, dt = \psi(q) + \gamma.
\end{equation}
\end{Example}

\begin{Example}
The evaluation $\mathbf{3.311.11}$:
\begin{equation}
\ift \frac{e^{px}-e^{qx}}{e^{rx}-e^{sx}} \, dx = 
\frac{1}{r-s} \left( \psi\left( \frac{r-q}{r-s} \right) - 
\psi \left( \frac{r-p}{r-s} \right) \right), 
\label{expquot}
\end{equation}
\noindent
follows directly from (\ref{psiaa}) by the change of variables $t=(r-s)x$. 
\end{Example}

\begin{Example}
The evaluation of $\mathbf{3.311.12}$:
\begin{equation}
\ift \frac{a^{x}-b^{x}}{c^{x}-d^{x}} \, dx = 
\frac{1}{\ln c - \ln d} 
\left( \psi\left( \frac{\ln c - \ln b}{\ln c - \ln d} \right) -
       \psi\left( \frac{\ln c - \ln a}{\ln c - \ln d} \right)  \right),
\end{equation}
\noindent
is proved by simply writing the exponentials in natural base. 
\end{Example}

\begin{Example}
The formula 
$\mathbf{3.311.10}$ had a sign error in the {\em sixth} edition of 
\cite{gr} : it appears as 
\begin{equation}
\ift \frac{e^{-px}-e^{-qx}}{1+e^{-(p+q)x}} \, dx = 
\frac{\pi}{p+q} \cot \left( \frac{p \pi}{p+q} \right). 
\end{equation}
\noindent
It should be 
\begin{equation}
\ift \frac{e^{-px}-e^{-qx}}{1-e^{-(p+q)x}} \, dx = 
\frac{\pi}{p+q} \cot \left( \frac{p \pi}{p+q} \right). 
\end{equation}
\noindent
The value (\ref{expquot}) yields
\begin{equation}
\ift \frac{e^{-px}-e^{-qx}}{1-e^{-(p+q)x}} \, dx =  
\frac{1}{p+q} \left( \psi\left( \frac{q}{p+q} \right) - 
 \psi\left( \frac{p}{p+q} \right) \right),
\end{equation}
\noindent
and the trigonometric answer follows from (\ref{psi2}). This has been 
corrected in the current edition of \cite{gr}.
\end{Example}

\begin{Example}
The evaluation of $\mathbf{3.312.2}$:
\begin{equation}
\ift \frac{(1-e^{-a x})(1- e^{-b x})e^{-px} }{1-e^{-x}} \, dx 
= \psi(p+a)+\psi(p+b) - \psi(p+a+b) - \psi(p)
\end{equation}
\noindent
follows directly from (\ref{3.231.5}). Indeed, the 
change of variables $t = e^{-x}$
gives
\begin{equation}
I = \int_{0}^{1} \frac{t^{p-1}(1-t^{a}-t^{b}+t^{a+b})}
{1-t} \, dt 
\end{equation}
\noindent
and now split them as
\begin{equation}
I = \int_{0}^{1} \frac{t^{p-1} - t^{p+a-1}}{1-t} \, dt -
    \int_{0}^{1} \frac{t^{p+b-1} - t^{p+a+b-1}}{1-t} \, dt
\end{equation}
\noindent
and use (\ref{3.231.5}) to conclude. 
\end{Example}

\section{A singular example} \label{sec-singular} 
\setcounter{equation}{0}

The example discussed in this section is
\begin{equation}
\int_{-\infty}^{\infty} \frac{e^{-\mu x} \, dx}{b - e^{-x}} 
= b^{\mu-1} \pi \, \cot( \pi \mu),
\end{equation}
\noindent
that appears as $\mathbf{3.311.8}$ in \cite{gr}. In the case $b > 0$ this has
to be modified in its presentation to avoid the
singularity $x = - \ln b$. The case $b < 0$ was discussed in 
\cite{moll-gr6}. In order to reduce the integral to a previous 
example, we let $t = e^{-x}$ to obtain 
\begin{equation}
\int_{-\infty}^{\infty} \frac{e^{-\mu x} \, dx}{b - e^{-x}} = 
\ift \frac{t^{\mu-1} \, dt}{b-t}.
\end{equation}
\noindent
The change of variables $t = by$ yields
\begin{equation}
\int_{-\infty}^{\infty} \frac{e^{-\mu x} \, dx}{b - e^{-x}} = 
b^{\mu-1} \ift \frac{y^{\mu-1} \, dy}{1-y}. 
\end{equation}
\noindent
Now separate the range of integration into $[0,1]$ and $[1,\infty)$. Then 
make the change of variables $y = 1/z$ in the second part. This produces 
\begin{equation}
\int_{-\infty}^{\infty} \frac{e^{-\mu x} \, dx}{b - e^{-x}} = 
b^{\mu-1} \int_{0}^{1} \frac{z^{\mu-1} - z^{-\mu}}{1 - z} \, dz. 
\end{equation}
\noindent
This last integral has been evaluated as $\cot(\pi \mu)$ in 
(\ref{3.231.1}).

\section{An integral with a fake parameter} \label{sec-fake1} 
\setcounter{equation}{0}

The example considered in this section is
$\mathbf{3.234.1}$:
\begin{equation}
\int_{0}^{1} \left( \frac{x^{q-1}}{1-ax} - \frac{x^{-q}}{a-x} \right) \, dx 
= \frac{\pi}{a^{q}} \cot \pi q.
\label{3.234.1}
\end{equation}
\noindent
We show that the parameter $a$ is {\em fake}, in the sense that it can 
be easily scaled out of the formula. The integral is written as 
$\lim\limits_{\epsilon \to 0} I(\epsilon)$ where
\begin{eqnarray}
I(\epsilon) & =  & 
\int_{0}^{1} \left( \frac{x^{q-1}}{(1-ax)^{1-\epsilon}} 
- \frac{x^{-q}}{(a-x)^{1- \epsilon}} \right) \, dx  \nonumber \\
& = & \int_{0}^{1} \frac{x^{q-1}}{(1-ax)^{1-\epsilon}} \, dx - 
\int_{0}^{1} \frac{x^{-q}}{(a-x)^{1- \epsilon}}  \, dx.  \nonumber 
\end{eqnarray}
\noindent
Make the change of variables $t =ax$ in the 
first integral and $x=at$ in the second one to produce 
\begin{equation}
I(\epsilon) = a^{-q} \int_{0}^{a} \frac{t^{q-1} \, dt}{(1-t)^{1- \epsilon}}-
a^{-q+\epsilon} \int_{0}^{1/a} \frac{t^{-q} \, dt}{(1-t)^{1 - \epsilon} },
\nonumber
\end{equation}
\noindent
and then let $\epsilon \to 0$ to produce
\begin{equation}
\int_{0}^{1} \left( \frac{x^{q-1}}{1-ax} - \frac{x^{-q}}{a-x} \right) \, dx  =
a^{-q} \left( \int_{0}^{a} \frac{t^{q-1} \, dt}{1-t} - 
\int_{0}^{1/a} \frac{t^{-q} \, dt}{1-t} \right). \nonumber
\end{equation}
\noindent
Differentiation with respect to the parameter $a$, shows that
the expression in parenthesis is independent of
$a$. It is now evaluated by using $a=1$ to obtain
\begin{equation}
\int_{0}^{1} \left( \frac{x^{q-1}}{1-ax} - \frac{x^{-q}}{a-x} \right) \, dx  =
a^{-q} \left( \int_{0}^{1} \frac{t^{q-1} -t^{-q}}{1-t} \, dt \right). \nonumber
\end{equation}
\noindent
The evaluation (\ref{3.231.5}) now yields 
\begin{eqnarray}
\int_{0}^{1} \left( \frac{x^{q-1}}{1-ax} - \frac{x^{-q}}{a-x} \right) \, dx  
& = & a^{-q} \left( \psi(1-q)-\psi(q) \right) \nonumber \\
& = & a^{-q} \pi \cot \pi q. \nonumber 
\end{eqnarray}
\noindent
Formula (\ref{3.234.1}) has been established. 

\medskip

\section{The derivative of  $\psi$} \label{sec-psiprime} 
\setcounter{equation}{0}

In a future publication we will discuss the evaluation of definite 
integrals in terms of the {\em polygamma function}
\begin{equation}
\text{PolyGamma}[n,x]:= \left( \frac{d}{dx} \right)^{n} \psi(x). 
\end{equation}
\noindent
In this section, we simply describe some integrals in \cite{gr} that comes 
from direct differentiation of the examples described above. 

\begin{Example}
Differentiating (\ref{3.231.5}) with respect to the parameter $p$ 
produces $\mathbf{4.251.4}$:
\begin{equation}
\int_{0}^{1} \frac{x^{p-1} \, \ln x }{1-x} \, dx = - \psi'(p).
\label{4.251.4}
\end{equation}
\end{Example}

\begin{Example}
The change of  variables $x = t^{q}$ in  (\ref{4.251.4}), followed by the
change of parameter $p \mapsto \frac{p}{q}$ yields $\mathbf{4.254.1}$:
\begin{equation}
\int_{0}^{1} \frac{t^{p-1} \, \ln t }{1-t^{q}} \, dx = - \frac{1}{q^{2}} 
\psi' \left( \frac{p}{q} \right).
\label{4.254.1}
\end{equation}
\end{Example}

\begin{Example}
Replace $q$ by $2q$ and $p$ by $q$ in (\ref{4.254.1}) to produce 
\begin{equation}
\int_{0}^{1} \frac{t^{q-1} \, \ln t }{1-t^{2q}} \, dt = - \frac{1}{4q^{2}}
\psi' \left( \tfrac{1}{2} \right). 
\end{equation}
\noindent
To evaluate this last term, differentiate the logarithm of the identity
\begin{equation}
\Gamma(2x) = \frac{2^{2x-1}}{\sqrt{\pi}} \Gamma(x) \, \Gamma(x+ \tfrac{1}{2}),
\end{equation}
\noindent
to obtain 
\begin{equation}
2 \psi(2x) = 2 \ln 2 + \psi(x) + \psi(x + \tfrac{1}{2} ). 
\end{equation}
\noindent
One more differentitation produces
\begin{equation}
4 \psi'(2x) = \psi'(x) + \psi'(x + \tfrac{1}{2} ).
\end{equation}
\noindent
The value $x = \tfrac{1}{2}$ gives 
\begin{equation}
\psi'( \tfrac{1}{2} ) =  3 \psi'(1) = \frac{\pi^{2}}{2}.
\end{equation}
\noindent
Therefore we obtain $\mathbf{4.254.6}$:
\begin{equation}
\int_{0}^{1} \frac{x^{q-1} \, \ln x }{1-x^{2q}} \, dx = - 
\frac{\pi^{2}}{8q^{2}}.
\end{equation}
\end{Example}

\begin{Example}
Differentiating (\ref{32443}) $n$-times with respect to  the parameter $p$ 
produces $\mathbf{4.271.15}$:
\begin{equation}
\int_{0}^{1} \ln^{n}x \, \frac{x^{p-1} \, dx}{1-x^{q}} = - 
\frac{1}{q^{n+1}} \psi^{(n)} \left( \frac{p}{q} \right).
\end{equation}
\end{Example}

\medskip

\section{A family of  logarithmic integrals} \label{sec-family} 
\setcounter{equation}{0}

Several of the integrals appearing in \cite{gr} are particular examples
of the family evaluated in the next proposition. 

\begin{Prop}
Let $a, \, b \in {\mathbb{R}}^{+}$. Then 
\begin{equation}
\int_{0}^{1} x^{a-1} (1-x)^{b-1} \, \ln x \, dx = 
\frac{\Gamma(a) \, \Gamma(b)}{\Gamma(a+b)} \left( \psi(a) - \psi(a+b) 
\right).
\label{form55}
\end{equation}
\end{Prop}
\begin{proof}
Differentiate the identity
\begin{equation}
\int_{0}^{1} x^{a-1} (1-x)^{b-1} \, dx  = 
\frac{\Gamma(a) \, \Gamma(b) }{\Gamma(a+b)}
\end{equation}
\noindent
with respect to the parameter $a$ and recall that 
$\Gamma'(x) = \psi(x) \Gamma(x)$. 
\end{proof}

\medskip

The next corollary appears as $\mathbf{4.253.1}$ in \cite{gr}. 

\begin{Cor}
\label{nice-coro}
Let $a, \, b, \, c \in \mathbb{R}^{+}$. Then
\begin{equation}
\int_{0}^{1} x^{a-1} (1-x^{c})^{b-1} \, \ln x \, dx = 
\frac{\Gamma(a/c) \, \Gamma(b)}{c^{2} \, \Gamma(a/c+b)} 
\left( \psi \left(\frac{a}{c} \right) - \psi \left(\frac{a}{c}+b \right) 
\right).
\label{form55a}
\end{equation}
\end{Cor}
\begin{proof}
Let $t = x^{c}$ in the integral (\ref{form55}).
\end{proof}

\begin{Example}
The formula in the previous corollary also appears as $\mathbf{4.256}$ in 
the form
\begin{equation}
\int_{0}^{1} \ln \left( \frac{1}{x} \right)  
\frac{x^{\mu-1} \, dx}{\sqrt[n]{(1-x^{n})^{n-m}}} = 
\frac{1}{n^{2}} B \left( \frac{\mu}{n}, \frac{m}{n} \right) 
\left[ \psi \left( \frac{\mu+m}{n} \right) - \psi \left( \frac{\mu}{n} \right) 
\right].
\end{equation}
\end{Example}

\begin{Example}
The integral 
\begin{equation}
\int_{0}^{1} \frac{x^{2n} \, \ln x}{\sqrt{1-x^{2}}} \, dx = 
\int_{0}^{1} x^{2n} (1-x^{2})^{-1/2} \, \ln x \, dx
\end{equation}
\noindent
that appears as $\mathbf{4.241.1}$ in \cite{gr}, corresponds to 
$a= 2n+1, \, b = \tfrac{1}{2}$ and $c=2$ in (\ref{form55a}). Therefore
\begin{equation}
\int_{0}^{1} \frac{x^{2n} \, \ln x}{\sqrt{1-x^{2}}} \, dx = 
\frac{\Gamma(n+\tfrac{1}{2}) \, \Gamma( \tfrac{1}{2} )}{4 \Gamma(n+1)} \, 
\left[ \psi(n + \tfrac{1}{2}) - \psi(n+1) \right].
\end{equation}
\noindent
Using (\ref{gamma-half}), (\ref{psi-n}) and (\ref{psi-half}) yields
\begin{equation}
\int_{0}^{1} \frac{x^{2n} \, \ln x}{\sqrt{1-x^{2}}} \, dx = 
\frac{\binom{2n}{n} \, \pi}{2^{2n+1}} 
\left( \sum_{k=1}^{2n} \frac{(-1)^{k-1}}{k} - \ln 2 \right). 
\label{42411}
\end{equation}
\noindent
This is $\mathbf{4.241.1}$.
\end{Example}

\begin{Example}
The integral in $\mathbf{4.241.2}$ states that
\begin{equation}
\int_{0}^{1} \frac{x^{2n+1} \, \ln x }{\sqrt{1-x^{2}}} \, dx = 
\frac{(2n)!!}{(2n+1)!!} 
\left( \ln 2 + \sum_{k=1}^{2n+1} \frac{(-1)^{k}}{k} \right). 
\label{42412old}
\end{equation}
\noindent
Writing the integral as
\begin{equation}
I = \int_{0}^{1} x^{2n+1} (1-x^{2})^{-1/2} \, \ln x \, dx 
\end{equation}
\noindent
we see that is corresponds to the case $a=2n+2, \, b = \tfrac{1}{2}, \, 
c= 2$ in (\ref{form55a}). Therefore
\begin{equation}
I = \frac{\Gamma(n+1) \, \Gamma( \tfrac{1}{2} ) }{4 \Gamma(n + \tfrac{3}{2})}
 \left[ \psi(n+1) - \psi(n + \tfrac{3}{2} ) \right].
\end{equation}
\noindent
Using (\ref{gamma-half}), (\ref{psi-n}) and (\ref{psi-half}) yields
\begin{equation}
\int_{0}^{1} \frac{x^{2n+1} \, \ln x }{\sqrt{1-x^{2}}} \, dx = 
\frac{2^{2n}}{(n+1) \, \binom{2n+1}{n} } 
\left( \ln 2 + \sum_{k=1}^{2n+1} \frac{(-1)^{k}}{k} \right)
\end{equation}
\noindent
This is equivalent to (\ref{42412old}).
\end{Example}

\begin{Example}
The integral $\mathbf{4.241.3}$ in \cite{gr} states that
\begin{equation}
\int_{0}^{1} x^{2n} \sqrt{1-x^{2}} \, \ln x \, dx = 
\frac{(2n-1)!!}{(2n+2)!!} \cdot \frac{\pi}{2} 
\left( \sum_{k=1}^{2n} \frac{(-1)^{k-1}}{k} - \frac{1}{2n+2} - \ln 2 \right).
\label{42413old}
\end{equation}
\noindent
To evaluate the integral, we write it as 
\begin{equation}
I = \int_{0}^{1} x^{2n} (1-x^{2})^{1/2} \, \ln x \, dx 
\end{equation}
\noindent
and we see that is corresponds to the case $a=2n+1, \, b = \tfrac{3}{2}, \, 
c= 2$ in (\ref{form55a}). Therefore
\begin{equation}
I = \frac{\Gamma(n+ \tfrac{1}{2} ) \, \Gamma(\tfrac{3}{2})}{4 \Gamma(n+2)}
\left[ \psi(n+ \tfrac{1}{2}) - \psi(n+2) \right].
\end{equation}
\noindent
Using (\ref{gamma-half}), (\ref{psi-n}) and (\ref{psi-half}) yields
\begin{equation}
\int_{0}^{1} x^{2n} \sqrt{1-x^{2}} \, \ln x \, dx = 
- \frac{\binom{2n}{n} \, \pi}{2^{2n+2} \, (n+1)} 
\left( \ln 2 + \frac{1}{2n+2} + \sum_{k=1}^{2n} \frac{(-1)^{k}}{k} \right).
\label{42413new}
\end{equation}
\noindent
This is equivalent to (\ref{42413old}).
\end{Example}

\begin{Example}
The integral $\mathbf{4.241.4}$ in \cite{gr} states that
\begin{equation}
\int_{0}^{1} x^{2n+1} \sqrt{1-x^{2}} \, \ln x \, dx = 
\frac{(2n)!!}{(2n+3)!!} 
\left( \ln 2 + \sum_{k=1}^{2n+1} \frac{(-1)^{k-1}}{k} - \frac{1}{2n+3} \right).
\label{42414old}
\end{equation}
\noindent
To evaluate the integral, we write it as 
\begin{equation}
I = \int_{0}^{1} x^{2n+1} (1-x^{2})^{1/2} \, \ln x \, dx 
\end{equation}
\noindent
and we see that is corresponds to the case $a=2n+2, \, b = \tfrac{3}{2}, \, 
c= 2$ in (\ref{form55a}). Therefore
\begin{equation}
I = \frac{\Gamma(n+ 1) \, \Gamma(\tfrac{3}{2})}{4 \Gamma(n+\tfrac{5}{2})}
\left[ \psi(n+ 1) - \psi(n+\tfrac{5}{2}) \right].
\end{equation}
\noindent
Using (\ref{gamma-half}), (\ref{psi-n}) and (\ref{psi-half}) yields
\begin{equation}
\int_{0}^{1} x^{2n+1} \sqrt{1-x^{2}} \, \ln x \, dx = 
\frac{2^{2n+1}}{(n+1)(n+2) \binom{2n+3}{n+1}}
\left( \ln 2 - \frac{1}{2n+3} + \sum_{k=1}^{2n+1} \frac{(-1)^{k}}{k} \right).
\label{42414new}
\end{equation}
\noindent
This is equivalent to (\ref{42414old}).
\end{Example}

\begin{Example}
The integral $\mathbf{4.241.5}$ in \cite{gr} states that
\begin{equation}
\int_{0}^{1} \ln x \, \sqrt{(1-x^{2})^{2n-1}}  \, dx = 
- \frac{(2n-1)!! \, \pi}{4 \, (2n)!!} \left[ \psi(n+1) + \gamma + \ln 4 
\right]
\label{42415old}
\end{equation}
\noindent
To evaluate the integral, we write it as 
\begin{equation}
I = \int_{0}^{1} (1-x^{2})^{n - \tfrac{1}{2}} \, \ln x \, dx 
\end{equation}
\noindent
and we see that is corresponds to the case $a=1, \, b = n+ \tfrac{1}{2}, \, 
c= 2$ in (\ref{form55a}). Therefore
\begin{equation}
I = \frac{\Gamma(n+ \tfrac{1}{2}) \, \Gamma(\tfrac{1}{2})}
{4 \Gamma(n+1)}
\left[ \psi(\tfrac{1}{2}) - \psi(n+\tfrac{1}{2}) \right].
\end{equation}
\noindent
Using (\ref{gamma-half}), (\ref{psi-n}) and (\ref{psi-half}) yields
\begin{equation}
\int_{0}^{1} (1-x^{2})^{n - \tfrac{1}{2}} \, \ln x \, dx 
= - \frac{\binom{2n}{n} \, \pi}{2^{2n+2}} 
\left( 2 \ln 2  + \sum_{k=1}^{n} \frac{1}{k} \right).
\end{equation}
\noindent
This is equivalent to (\ref{42415old}). This integral also appears as 
$\mathbf{4.246}$. 
\end{Example}

\begin{Example}
The case $n=0$ in (\ref{42411}) yields
\begin{equation}
\int_{0}^{1} \frac{\ln x \, dx }{\sqrt{1-x^{2}}} = - \frac{\pi}{2} \ln 2.
\end{equation}
\noindent
This appears as $\mathbf{4.241.7}$ in \cite{gr}.
\end{Example}

\begin{Example}
Formula $\mathbf{4.241.8}$ states that
\begin{equation}
\int_{1}^{\infty} \frac{\ln x \, dx}{x^{2} \, \sqrt{x^{2}-1} } = 1 - \ln 2.
\end{equation}
\noindent
To evaluate this, let $t = 1/x$ to obtain
\begin{equation}
 I = - \int_{0}^{1} t (1-t^{2})^{-1/2} \, \ln t  \, dt.
\end{equation}
\noindent
This corresponds to the case $a=2, \, b = \tfrac{1}{2}, \, c=2$ in 
(\ref{form55a}). Therefore
\begin{equation}
I = - \frac{\Gamma(1) \, \Gamma( \tfrac{1}{2} ) }{4 \, \Gamma(\tfrac{3}{2})}
\left[ \psi(1) - \psi( \tfrac{3}{2} ) \right],
\end{equation}
\noindent
and the value $1- \ln 2$ comes from (\ref{psi-n}) and (\ref{psi-half}).
\end{Example}

\begin{Example}
The case $n=0$ in (\ref{42413new}) produces
\begin{equation}
\int_{0}^{1} \sqrt{1-x^{2}} \, \ln x \, dx = - \frac{\pi}{8}(2 \ln 2 + 1).
\end{equation}
\noindent
This appears as $\mathbf{4.241.9}$ in \cite{gr}.
\end{Example}

\begin{Example}
The case $n=0$ in (\ref{42414new}) produces
\begin{equation}
\int_{0}^{1} x \, \sqrt{1-x^{2}} \, \ln x \, dx =  \frac{1}{9}(3 \ln 2 - 4).
\end{equation}
\noindent
This appears as $\mathbf{4.241.10}$ in \cite{gr}.
\end{Example}

\begin{Example}
Entry $\mathbf{4.241.11}$ states that
\begin{equation}
\int_{0}^{1} \frac{ \ln x \, dx}{\sqrt{x(1-x^{2})}} = - 
\frac{\sqrt{2 \pi}}{8} \Gamma^{2} \left( \tfrac{1}{4} \right).
\end{equation}
\noindent
To evaluate the integral, write it as 
\begin{equation}
I = \int_{0}^{1} x^{-1/2} (1-x^{2})^{-1/2} \, \ln x \, dx 
\end{equation}
\noindent
and this corresponds to the case $a= \tfrac{1}{2}, \, b = \tfrac{1}{2}, \, 
c= 2$ in (\ref{form55a}). Therefore
\begin{equation}
I = \frac{\Gamma( \tfrac{1}{4} ) \, \Gamma( \tfrac{1}{2} )}
{4 \Gamma( \tfrac{3}{4} )} \left[ \psi \left( \tfrac{1}{4} \right)  - 
\psi \left( \tfrac{3}{4} \right) \right].
\end{equation}
\noindent
The stated form comes from using (\ref{reflec-gamma}) and (\ref{psi2}).
\end{Example}

\begin{Example}
The identity
\begin{equation}
\int_{0}^{1} \frac{x \, \ln x }{\sqrt{1-x^{4}}} \, dx = - \frac{\pi}{8} \ln 2
\end{equation}
\noindent
appears as $\mathbf{4.243}$ in \cite{gr}. To evaluate it, we write it as
\begin{equation}
I = \int_{0}^{1} x (1-x^{4})^{-1/2} \, \ln x \, dx 
\end{equation}
\noindent
that corresponds to $a =2, \, b = \tfrac{1}{2}, \, c=4$ in (\ref{form55a}).
Therefore,
\begin{equation}
I = \frac{1}{16} \Gamma^{2}( \tfrac{1}{2} ) \left[ \psi(\tfrac{1}{2}) - 
\psi(1) \right]. 
\end{equation}
\noindent
The values $\psi(1) = - \gamma$ and $\psi( \tfrac{1}{2} ) = - \gamma - 2 \ln 2$
gives the result.
\end{Example}

\begin{Example}
The verification of $\mathbf{4.244.1}$:
\begin{equation}
\int_{0}^{1} \frac{\ln x \, dx}{\sqrt[3]{x(1-x^{2})^{2} } } = 
- \frac{1}{8} \Gamma^{3} \left( \tfrac{1}{3} \right)
\end{equation}
\noindent
is achieved by using (\ref{form55a}) with $a= \tfrac{2}{3}, \, b= \tfrac{1}{3}$
and $c=2$ to obtain
\begin{equation}
I = \frac{\Gamma^{2} \left( \tfrac{1}{3} \right)}{4 \Gamma( \tfrac{2}{3} )} 
\left[ \psi \left( \tfrac{1}{3} \right)  - \psi \left(\tfrac{2}{3} \right)
\right].
\end{equation}
\noindent
Using (\ref{reflec-gamma}) and (\ref{psi2}) produces the stated result.
\end{Example}

\begin{Example}
The usual application of (\ref{form55a}) shows that $\mathbf{4.244.2}$ is
\begin{equation}
\int_{0}^{1} \frac{\ln x \, dx}{\sqrt[3]{1-x^{3}}} = 
\frac{2 \pi}{9 \sqrt{3}} \left( \psi \left( \tfrac{1}{3} \right) + 
\gamma \right),
\end{equation}
\noindent
where we have used 
$\Gamma( \tfrac{1}{3} ) \Gamma( \tfrac{2}{3}) = 2 \pi/\sqrt{3}$. It remains 
to evaluate $\psi(\tfrac{1}{3})$.  The identity (\ref{psi2})  gives
\begin{equation}
\psi \left( \tfrac{1}{3} \right) - \psi \left( \tfrac{2}{3} \right) = 
- \frac{\pi}{\sqrt{3}}.
\end{equation}
\noindent
To obtain a second relation among these quantities, we start from the 
identity 
\begin{equation}
\Gamma(3x) = \frac{3^{3x-1/2}}{2 \pi} 
\Gamma(x) \, \Gamma(x + \tfrac{1}{3}) \, \Gamma( x + \tfrac{2}{3} )
\end{equation}
\noindent
that follows directly from  (\ref{gamma-product}), and differentiate 
logarithmically to obtain
\begin{equation}
\psi(3x) = \ln 3 + \frac{1}{3} \left( \psi(x) + \psi( x + \tfrac{1}{3} ) +
\psi( x + \tfrac{2}{3} \right). 
\end{equation}
\noindent
The special case $x = \tfrac{1}{3}$ yields
\begin{equation}
\psi \left( \tfrac{1}{3} \right) + \psi \left( \tfrac{2}{3} \right) = 
-2 \gamma - 3 \ln 3. 
\end{equation}
\noindent
We conclude that
\begin{equation}
\psi( \tfrac{1}{3} ) = - \gamma - \frac{3}{2} \ln 3 - \frac{\pi}{2 \sqrt{3}}
\label{psi13}
\end{equation}
\noindent
and
\begin{equation}
\psi( \tfrac{2}{3} ) = - \gamma - \frac{3}{2} \ln 3 + \frac{\pi}{2 \sqrt{3}}.
\label{psi23}
\end{equation}
\noindent
This gives
\begin{equation}
\int_{0}^{1} \frac{\ln x \, dx}{\sqrt[3]{1-x^{3}}} = 
- \frac{\pi}{3 \sqrt{3}} \left( \ln 3 + \frac{\pi}{3 \sqrt{3}} \right),
\end{equation}
\noindent
as stated in $\mathbf{4.244.2}$. 
\end{Example}

\begin{Example}
The evaluation of $\mathbf{4.244.3}$:
\begin{equation}
\int_{0}^{1} \frac{\ln x \, dx}{\sqrt[3]{1-x^{3}}} = 
- \frac{\pi}{3 \sqrt{3}} \left( \ln 3 - \frac{\pi}{3 \sqrt{3}} \right),
\end{equation}
\noindent
proceeds as in the previous example. The integral is identified as 
\begin{equation}
I = \tfrac{1}{9} \Gamma( \tfrac{2}{3} ) \Gamma( \tfrac{1}{3} ) 
\left[ \psi \left( \tfrac{2}{3} \right) + \gamma \right].
\end{equation}
\noindent
The value (\ref{psi23}) gives the rest.
\end{Example}

\begin{Example}
The change of variables $t=x^{4}$ yields
\begin{equation}
\int_{0}^{1} \frac{x^{p} \, \ln x \, dx }{\sqrt{1-x^{4}}} = 
\frac{1}{16} \int_{0}^{1} t^{(p-3)/4} (1-t)^{-1/2} \, \ln t \, dt.
\end{equation}
\noindent
The last integral is evaluated using (\ref{form55a}) with $a = \tfrac{p+1}{4},
\, b = \tfrac{1}{2}$ and $c=1$ to obtain
\begin{equation}
\int_{0}^{1} \frac{x^{p} \, \ln x \, dx }{\sqrt{1-x^{4}}} = 
\frac{\sqrt{\pi}}{16} \frac{\Gamma( \tfrac{p+1}{4})}{\Gamma( \tfrac{p+3}{4})}
\left[ \psi \left( \frac{p+1}{4} \right) - \psi \left( \frac{p+3}{4} \right)
\right].
\end{equation}
\noindent
The special case $p= 4n+1$ yields
\begin{equation}
\int_{0}^{1} \frac{x^{4n+1} \, \ln x \, dx }{\sqrt{1-x^{4}}} = 
\frac{\sqrt{\pi}}{16n!} \Gamma( n + \tfrac{1}{2} ) \left[ 
\psi(n+ \tfrac{1}{2}) - \psi(n+1) \right].
\end{equation}
The special case $p= 4n+1$ yields
\begin{equation}
\int_{0}^{1} \frac{x^{4n+1} \, \ln x \, dx }{\sqrt{1-x^{4}}} = 
\frac{\sqrt{\pi} \, \Gamma(n+ \tfrac{1}{2}) }{16 n!} 
\left[ \psi(n+ \tfrac{1}{2} ) - \psi(n+1) \right].
\nonumber
\end{equation}
\noindent
Using (\ref{gamma-half}), (\ref{psi-n}) and (\ref{psi-half}) yields
$\mathbf{4.245.1}$ in the form
\begin{equation}
\int_{0}^{1} \frac{x^{4n+1} \, \ln x \, dx }{\sqrt{1-x^{4}}} = 
\frac{\pi \, \binom{2n}{n} }{2^{2n+3}} \left( - \ln 2 + 
\sum_{k=1}^{2n} \frac{(-1)^{k-1}}{k} \right).
\end{equation}
\noindent
The special case $p= 4n+3$ yields
\begin{equation}
\int_{0}^{1} \frac{x^{4n+3} \, \ln x \, dx }{\sqrt{1-x^{4}}} = 
\frac{\sqrt{\pi} \, n!}{16 \Gamma( n + \tfrac{3}{2})} 
\left[ \psi( n+1) - \psi(n + \tfrac{3}{2} ) \right].
\nonumber
\end{equation}
Using (\ref{gamma-half}), (\ref{psi-n}) and (\ref{psi-half}) yields
$\mathbf{4.245.2}$ in the form
\begin{equation}
\int_{0}^{1} \frac{x^{4n+3} \, \ln x \, dx }{\sqrt{1-x^{4}}} = 
\frac{2^{2n-2}}{(2n+1) \, \binom{2n}{n} } 
\left( \ln 2 + \sum_{k=1}^{2n+1} \frac{(-1)^{k}}{k} \right).
\end{equation}
\end{Example}

\begin{Example}
The change of variables $t = x^{2n}$ produces
\begin{equation}
\int_{0}^{1} \frac{\ln x \, dx }{\sqrt[n]{1-x^{2n}}} = 
\frac{1}{4n^{2}} \int_{0}^{1} t^{\tfrac{1}{2n} -1} (1-t)^{-\tfrac{1}{n}} 
\, \ln t \, dt.
\end{equation}
\noindent
Then (\ref{form55a}) with $a = \tfrac{1}{2n}, \, b = 1 - \tfrac{1}{n}$ 
and $c= 1$ give the value
\begin{equation}
\int_{0}^{1} \frac{\ln x \, dx }{\sqrt[n]{1-x^{2n}}} = 
\frac{\Gamma(\tfrac{1}{2n} ) \, \Gamma(1  - \tfrac{1}{n} )}
{4n^{2} \, \Gamma(1 - \tfrac{1}{2n}) }
\left[ \psi \left( \tfrac{1}{2n} \right) - \psi \left( 1 - \tfrac{1}{2n} 
\right) \right]. 
\end{equation}
\noindent
Using (\ref{gamma-half}) and (\ref{psi-half}) to obtain
\begin{equation}
\int_{0}^{1} \frac{\ln x \, dx}{\sqrt[n]{1-x^{2n}}} = 
-\frac{\pi}{8} \frac{B \left( \tfrac{1}{2n}, \tfrac{1}{2n} \right)}
{n^{2} \, \sin \left( \tfrac{\pi}{2n} \right)}.
\end{equation}
\noindent
This is $\mathbf{4.247.1}$ in \cite{gr}.
\end{Example}

\begin{Example}
The change of variables $t = x^{2}$ gives
\begin{equation}
\int_{0}^{1} \frac{\ln x \, dx}{\sqrt[n]{x^{n-1} (1-x^{2})}} = 
\frac{1}{4} \int_{0}^{1} t^{\tfrac{1}{2n}-1} (1-t)^{- \tfrac{1}{n}} \, 
\ln t \, dt. \nonumber
\end{equation}
\noindent
Using (\ref{form55a}) we obtain
\begin{equation}
\int_{0}^{1} \frac{\ln x \, dx}{\sqrt[n]{x^{n-1} (1-x^{2})}} = 
\frac{\Gamma( \tfrac{1}{2n}) \, \Gamma( 1 - \tfrac{1}{2n} )}
{ \Gamma( 1 - \tfrac{1}{2n} )} 
\left[ \psi \left( \tfrac{1}{2n} \right) - 
       \psi \left( 1- \tfrac{1}{2n} \right)  \right].
\end{equation}
\noindent
Proceeding as in the previous example, we obtain
\begin{equation}
\int_{0}^{1} \frac{\ln x \, dx}{\sqrt[n]{x^{n-1} (1-x^{2})}} = 
-\frac{\pi}{8} \frac{B \left( \tfrac{1}{2n}, \tfrac{1}{2n} \right)}
{ \sin \left( \tfrac{\pi}{2n} \right)}.
\end{equation}
\noindent
This is $\mathbf{4.247.2}$ in \cite{gr}.
\end{Example}

Some integrals in \cite{gr} have the form  of the Corollary \ref{nice-coro}
after an elementary change of variables. 

\begin{Example}
Formula $\mathbf{4.293.8}$ in \cite{gr} states that 
\begin{equation}
\int_{0}^{1} x^{a-1} \ln(1-x) \, dx  = -\frac{1}{a} \left( \psi(a+1) + 
\gamma \right).
\end{equation}
This follows directly from (\ref{form55a}) by the change of variables 
$x \mapsto 1-x$.  The same is true for $\mathbf{4.293.13}$:
\begin{equation}
\int_{0}^{1} x^{a-1} (1-x)^{b-1} \, \ln (1-x) \, dx = 
B(a,b) \left[ \psi(b) - \psi(a+b) \right].
\end{equation}
\end{Example}

\begin{Example}
The  change of variables $t = e^{-x}$ gives
\begin{equation}
\int_{0}^{\infty} xe^{-x} ( 1 - e^{2x})^{n - \tfrac{1}{2}} \, dx = 
- \int_{0}^{1} (1-t^{2})^{n - \tfrac{1}{2}} \, \ln t \, dt.
\end{equation}
\noindent
This latter integral is evaluated using (\ref{form55a}) as 
\begin{equation}
I = - \frac{\sqrt{\pi} \Gamma( n + \tfrac{1}{2} )}{4 n!} 
\left( \psi\left( \tfrac{1}{2} \right) - \psi(n+1) \right).
\end{equation}
\noindent
Using (\ref{gamma-half}) and (\ref{psi-n}) we obtain 
\begin{equation}
\int_{0}^{\infty} xe^{-x} ( 1 - e^{-2x})^{n - \tfrac{1}{2}} \, dx = 
\frac{\binom{2n}{n} \, \pi}{2^{2n+2}} \left( 2 \ln 2 + \sum_{k=1}^{n} 
\frac{1}{k} \right). 
\end{equation}
\noindent
This appears as $\mathbf{3.457.1}$ in \cite{gr}. 
\end{Example}

\section{An announcement} \label{sec-announcement} 
\setcounter{equation}{0}

There are many integrals in \cite{gr} that contain the 
term $1+x$ in the denominator, instead of the term $1-x$ seen, for instance, 
in Section \ref{sec-simple}. The evaluation of these integrals can be 
obtained using the {\em incomplete beta function}, defined by
\begin{equation}
\beta(x) := \int_{0}^{1} \frac{t^{x-1} \, dt}{1+x}
\end{equation}
\noindent
as it appears in $\mathbf{8.371.2}$. This function is related to the 
digamma function by the identity
\begin{equation}
\beta(x) = \frac{1}{2} \left[ \psi \left( \frac{x+1}{2} \right) - 
\psi \left( \frac{x}{2} \right) \right]. 
\end{equation}
\noindent
These evaluations will be reported in \cite{moll-gr11}.

\section{One more family} \label{sec-onemore} 
\setcounter{equation}{0}

We conclude this collection with a two-parameter family of  integrals.

\begin{Prop}
Let $a, \, b \in \mathbb{R}^{+}$. Then 
\begin{equation}
\int_{0}^{\infty} \left( e^{-x^{a}} - \frac{1}{1+x^{b}} \right) \, \frac{dx}{x}
= - \frac{\gamma}{a},
\end{equation}
\noindent
independently of $b$.
\end{Prop}
\begin{proof}
Write the integral as 
\begin{equation}
\int_{0}^{\infty} \left( e^{-x^{a}} - e^{-x^{b}} \right) \, \frac{dx}{x} + 
\int_{0}^{\infty} \left( e^{-x^{b}} - \frac{1}{1+x^{b}}  \right) \, 
\frac{dx}{x}. 
\end{equation}
\noindent
The first integral is $(a-b)\gamma/ab$ according to (\ref{34762}). The change
of variables $t = x^{b}$ converts the second one into
\begin{equation}
\frac{1}{b} \int_{0}^{\infty} \left( e^{-t} - \frac{1}{1+t} \right) 
\, \frac{dt}{t}  = - \frac{\gamma}{b},
\end{equation}
\noindent
according to (\ref{34353}). The formula has been established.
\end{proof}

\begin{Example}
The  case $a= 2^{n}$ and $b = 2^{n+1}$ gives $\mathbf{3.475.1}$:
\begin{equation}
\int_{0}^{\infty} \left( \text{exp}(-x^{2^{n}}) - \frac{1}{1+x^{2^{n+1}}} 
\right) \, \frac{dx}{x} = -\frac{\gamma}{2^{n}}.
\end{equation}
\end{Example}

\begin{Example}
The case $a=2^{n}$ and $b=2$ gives $\mathbf{3.475.2}$:
\begin{equation}
\int_{0}^{\infty} \left( \text{exp}(-x^{2^{n}}) - \frac{1}{1+x^{2}} 
\right) \, \frac{dx}{x} = -\frac{\gamma}{2^{n}}.
\end{equation}
\end{Example}

\begin{Example}
The case $a=2$ and $b=2$ gives $\mathbf{3.467}$:
\begin{equation}
\int_{0}^{\infty} \left( e^{-x^{2}}  - \frac{1}{1+x^{2}} 
\right) \, \frac{dx}{x} = -\frac{\gamma}{2}.
\end{equation}
\end{Example}

\begin{Example}
Finally, the change of variables $t = ax$ yields 
\begin{equation}
\int_{0}^{\infty} \left( e^{-px} - \frac{1}{1+a^{2}x^{2} } \right) \frac{dx}{x} 
 = 
\int_{0}^{\infty} \frac{e^{-pt/a} - e^{-t}}{t} \, dt +
\int_{0}^{\infty}  \left( e^{-t} - \frac{1}{1+t^{2}} \right) \frac{dt}{t}.
\end{equation}
\noindent
The first integral is $\ln \frac{a}{p}$ according to (\ref{form25}), the 
second one is $-\gamma$. This gives the evaluation of $\mathbf{3.442.3}$:
\begin{equation}
\int_{0}^{\infty} \left( e^{-px} - \frac{1}{1+a^{2}x^{2} } \right) \frac{dx}{x} 
 = \gamma + \ln \frac{a}{p}.
\end{equation}
\end{Example}


\bigskip
\begin{thebibliography}{1}

\bibitem{irrbook}
G.~Boros and V.~Moll.
\newblock {\em Irresistible {I}ntegrals}.
\newblock Cambridge {U}niversity {P}ress, {N}ew {Y}ork, 1st edition, 2004.

\bibitem{gr}
I.~S. Gradshteyn and I.~M. Ryzhik.
\newblock {\em Table of {I}ntegrals, {S}eries, and {P}roducts}.
\newblock Edited by A. Jeffrey and D. Zwillinger. Academic Press, New York, 7th
  edition, 2007.

\bibitem{moll-gr4}
V.~Moll.
\newblock The integrals in {G}radshteyn and {R}yzhik. {P}art 4: {T}he gamma
  function.
\newblock {\em Scientia}, 15, 2007.

\bibitem{moll-gr6}
V.~Moll.
\newblock The integrals in {G}radshteyn and {R}yzhik. {P}art 6: {T}he beta
  function.
\newblock {\em Scientia}, to appear.

\bibitem{moll-gr11}
V.~Moll.
\newblock The integrals in {G}radshteyn and {R}yzhik. {P}art 11: {T}he
  incomplete beta function.
\newblock {\em Scientia}, in preparation.

\end{thebibliography}
\end{document}